\documentclass{article}
\usepackage{latexsym}
\usepackage{epic,eepic,color}
\usepackage{graphicx}
\usepackage{graphics}
\usepackage{amsmath}
\usepackage{amssymb}

\def\rcr{\overline{\hbox{\rm cr}}}

\def\enetercios{{s}}

\def\digaawithpi{{D_\Pi^{aa}}}
\def\digabwithpi{{D_\Pi^{bb}}}
\def\digacwithpi{{D_\Pi^{cc}}}

\def\digraa{{D^{aa}}}
\def\digab{{D^{bb}}}
\def\digac{{D^{cc}}}

\def\dd{{\cal D}}

\def\nomberle#1#2{{\chi_{\le #1}(#2)}}
\def\Nomberle#1#2{{\eta_{\le #1}(#2)}}
\def\Nomberlehom#1#2{{\eta_{\le #1}^{{}^{\hbox{\small\rm  hom}}}(#2)}}
\def\Nomberlehet#1#2{{\eta_{\le #1}^{{}^{\hbox{\small\rm het}}}(#2)}}
\def\Momberaa#1#2{{\eta_{> #1}^{{{aa}}}(#2)}}
\def\Momberbb#1#2{{\eta_{> #1}^{{{bb}}}(#2)}}
\def\Mombercc#1#2{{\eta_{> #1}^{{{cc}}}(#2)}}

\def\theb{{\rm b}}
\def\ther{{\rm r}}
\def\theq{{\rm q}}

\def\ygriega#1#2{{Y(#1,#2)}}

\def\outd#1{{[#1]^+}}
\def\outdpi#1{{[#1]_\Pi^+}}
\def\ind#1{{[#1]^-}}
\def\indpi#1{{[#1]_{\Pi}^-}}

\def\ora#1{{\overrightarrow{#1}}}

\def\hal{\vrule height 5 pt width 0.05 in  depth 0.8 pt}

\newtheorem{theorem}{Theorem}
\newtheorem{corollary}[theorem]{Corollary}

\newtheorem{claim}[theorem]{Claim}
\newtheorem{proposition}[theorem]{Proposition}

\newtheorem{remark}[theorem]{Remark}

\begin{document}

\author{B.~\'Abrego\footnote{Department of
      Mathematics, California State University Northridge.} %
\and S.~Fern\'andez--Merchant$^*$ 
\and J.~Lea\~nos\footnote{Instituto de F\'\i sica, Universidad
    Aut\'onoma de San Luis Potos\'{\i}, Mexico. Supported by
    FAI--UASLP and by CONACYT Grant 45903.}
\and G.~Salazar$^\dagger$}

\title{On $3$--decomposable geometric drawings of $K_n$}

\maketitle


\begin{abstract}
The point sets of all known optimal rectilinear drawings of $K_n$
share an unmistakeable 
clustering property, the so--called {\em 3--decomposability}.
It is widely believed that the underlying point sets of 
all optimal rectilinear drawings of $K_n$
are $3$--decomposable.
We give a lower bound for the minimum number of $(\le k)$--sets in a
$3$--decomposable $n$--point set. As an immediate corollary, we obtain a lower
bound for the crossing number $\rcr(\dd)$ of any
rectilinear drawing $\dd$ of $K_n$ with underlying
$3$--decomposable point set, namely
$\rcr(\dd) > \frac{2}{27}\left(15-\pi^{2}\right)\binom{n}{4}+\Theta(n^{3})
\approx 0.380029\binom{n}{4} + \Theta(n^3)$.
This closes this gap between the best known lower and upper bounds for
the rectilinear crossing number $\rcr(K_n)$ of $K_n$ 
by over 40\%, under the  assumption of
$3$--decomposability. 
\end{abstract}

\section{Introduction}

Figure~\ref{fig:figure1} shows the point set of an optimal 
(crossing minimal) rectilinear drawing of
$K_9$, with an  evident partition of the $9$ vertices into $3$
highly structured clusters of $3$ vertices each:

\begin{figure}[ht]\label{fig:figure1}
\begin{center}
\includegraphics[width=3cm]{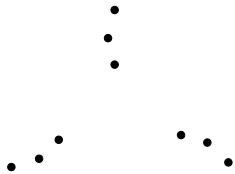}
\end{center}
\caption{The points in this optimal drawing of $K_9$ are
clustered into $3$ sets.}
\end{figure}

A similar, natural, highly structured partition into $3$ clusters of equal size is
observed in {\em every} known optimal drawing of $K_n$, for every
$n$  multiple of $3$~(see \cite{aweb}).
Even for those values of $n$  (namely,  $n > 27$) for which the
exact rectilinear crossing number $\rcr(K_n)$ of $K_n$ is not known,
the best available examples also share this property~\cite{aweb}.

In all these examples, a set $S$ of $n$
points in general position is partitioned into sets $A, B,$ and $C$,
with $|A| = |B| = |C| = n/3$ with the following properties:

\begin{description}
\item{(i)} There is a directed line $\ell_1$ such that, as we traverse
$\ell_1$, we find  the $\ell_1$--orthogonal
projections of the points in $A$, then the $\ell_1$--orthogonal
projections of the points in $B$, and then the $\ell_1$--orthogonal
projections of the points in $C$;
\item{(ii)} there is a directed line $\ell_2$ such that, as we traverse
$\ell_2$, we find  the $\ell_2$--orthogonal
projections of the points in $B$, then the $\ell_2$--orthogonal
projections of the points in $A$, and then the $\ell_2$--orthogonal
projections of the points in $C$; and
\item{(iii)} there is a directed line $\ell_3$ such that, as we traverse
$\ell_3$, we find  the $\ell_3$--orthogonal
projections of the points in $B$, then the $\ell_3$--orthogonal
projections of the points in $C$, and then the $\ell_3$--orthogonal
projections of the points in $A$.
\end{description}

\noindent{\bf Definition }
A point set that satisfies conditions (i)--(iii) above is $3$--{\em
  decomposable}.  We also say that the underlying rectilinear drawing
  of $K_n$ is $3$--{\em decomposable}.

\bigskip

\begin{figure} [ht]
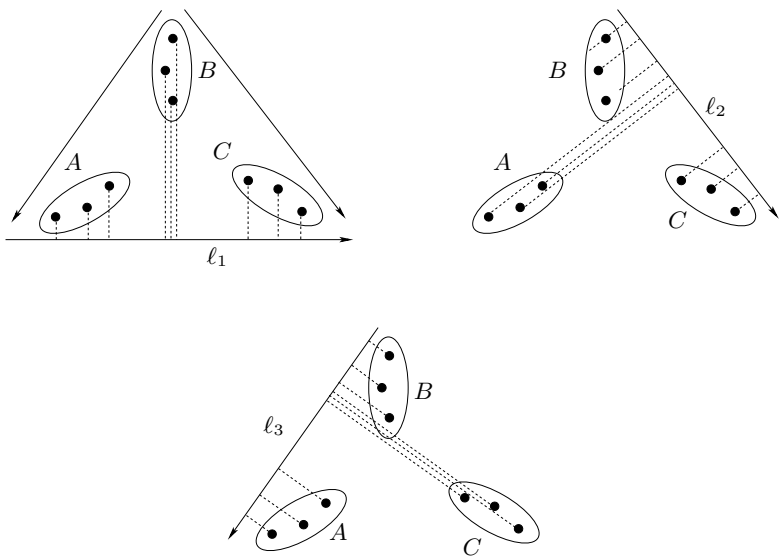

\centering
\input{k9-3l1.pstex_t} \hglue 1 cm
\input{k9-3l2.pstex_t} \vglue 0.5 cm
\input{k9-3l3.pstex_t}
\caption{The $9$--point set $S$ gets naturally partitioned into three
  clusters $A, B,$ and $C$. The $\ell_1$--, $\ell_2$--, and
  $\ell_3$--orthogonal projections of $A, B$, and $C$ satisfy
  conditions (i)--(iii), and so $S=A\cup B \cup C$ is $3$--decomposable.}
\label{fig:figure2}
\end{figure}

A possible choice of $\ell_1, \ell_2$, and $\ell_3$ for the example
in Figure~\ref{fig:figure1} is illustrated in Figure~\ref{fig:figure2}.

\subsection{The main result}

It is widely believed that all optimal rectilinear drawings of $K_n$
are $3$--decomposable.  One of our main results in this paper is the
following lower bound for the number of crossings in all such
drawings.

\begin{theorem}\label{thm:main}
Let $\dd$ be a $3$--decomposable rectilinear drawing of $K_n$. Then
the number $\rcr(\dd)$ of crossings in $\dd$ satisfies
$$\rcr(\dd) \ge 
\frac{2}{27}\left(15-\pi^{2}\right)\binom{n}{4}+\Theta(n^{3})
\ \approx \  0.380029\binom{n}{4} + \Theta(n^3).$$
\end{theorem}

The best known general lower and upper bounds for the rectilinear
crossing number $\rcr(K_n)$ are $ 0.37968\binom{n}{4} + \Theta(n^3)
\le \rcr(K_n) \le 0.38054\binom{n}{4} + \Theta(n^3)$ (see~\cite{agor} 
and~\cite{seq}).  Thus the bound given by Theorem~\ref{thm:main}
closes this gap by over 40\%, under the (quite feasible) assumption of
$3$--decomposability. 

To prove Theorem~\ref{thm:main} (in Section~\ref{sec:proofmain}), 
we exploit the close relationship
between rectilinear crossing numbers and $(\le k)$--sets,
unveiled independently by \'Abrego and Fern\'andez--Merchant~\cite{af}
and by Lov\'asz et al.~\cite{lovasz}.  

Recall that a $(\le k)$--{\em set} of a point set $S$ is a subset $T$
of $S$ with $|T| \le k$ such that some straight line separates $T$ and
$S\setminus T$.  The number $\nomberle{k}{S}$ of $(\le k)$--sets of
$S$ is a parameter of independent interest in discrete geometry
(see~\cite{bmp}), and, as we recall in Section~\ref{sec:proofmain}, is
closely related to the rectilinear crossing number of the geometric
graph induced by $S$.

The main ingredient in the proof of Theorem~\ref{thm:main} is the
following bound (Theorem~\ref{thm:mainksets}) for the 
number of $(\le k)$--sets in $3$--decomposable
point sets. The bound is in terms of the following quantity (by
convention, $\binom{r}{s} = 0$ 
if $r < s$), 

\begin{equation}\label{eq:ygriega}
\ygriega{k}{n} := 3\binom{k+1}{2}+3\binom{k-n/3+1}{2}+3\sum_{j=2}^{\theb}
j(j+1)\binom{k+1-(\frac{1}{2}-\frac{1}{3j(j+1)})3n}{2} 
- \frac{1}{3},
\end{equation}
where  $\theb:=\theb(k,n)$ is the unique integer such that 
$\binom{\theb(k,n) + 1}{2} <  n/(n-2k-1)   \le \binom{\theb(k,n) +2}{2}$.

\bigskip

\begin{theorem}\label{thm:mainksets}
Let $S$ be a $3$--decomposable set of 
$n$ points in general position, where $n$ is a multiple of $3$, and
let $k < n/2$.    
Then
$$
\nomberle{k}{S}  \ge \ygriega{k}{n}.
$$
\end{theorem}

\bigskip

The best general lower bound for $\nomberle{k}{S}$ is
the sum of the first two terms in (\ref{eq:ygriega}) (see~\cite{agor}
and~\cite{seq}). Thus the third
summand in (\ref{eq:ygriega}) is the improvement we report, under the
assumption of $3$--decomposability.

The proofs of Theorems~\ref{thm:main} and~\ref{thm:mainksets} are in
Sections~\ref{sec:proofmain} and~\ref{sec:proofmainksets},
respectively. In Section~\ref{sec:concludingremarks} we present some
concluding remarks and open questions.

\section{Proof of Theorem~\ref{thm:main}}\label{sec:proofmain}

Let $\dd$ be a $3$--decomposable rectilinear drawing of $K_n$, and
let $S$ denote the underlying $n$--point set, that is, the vertex set
of $\dd$.  Besides Theorem~\ref{thm:mainksets}, our main tool is the
following relationship between $(\le k)$--sets and the rectilinear
crossing number (see~\cite{af} or~\cite{lovasz}):
\begin{equation}\label{eq:aflov}
\rcr(\dd)=%
\sum_{1\leq k\leq(n-2)/2}
(n-2k-1)\nomberle{k}{S}+\Theta(n^{3}).%
\end{equation}

Combining Theorem~\ref{thm:mainksets} and Eq.~(\ref{eq:aflov}), and
noting that both the $-1$ in the factor $n-2k-1$ and the summand
$-1/3$ in (\ref{eq:ygriega})
only contribute to smaller order
terms, we obtain:

\begin{align*}
&\textstyle{\rcr(\dd_n)}
\textstyle{\ge\sum_{k=1}^{(n-2)/2} (n-2k)  
\biggl(  
3\binom{k+1}{2}+3\binom{k-n/3+1}{2}+3\sum_{j=2}^{\theb}
j(j+1)\binom{k+1-(\frac{1}{2}-\frac{1}{3j(j+1)})3n}{2} 
\biggr)
+\Theta(n^3)} \\
& \textstyle{\ge 24{n\choose 4}} \sum_{k=1}^{(n-2)/2} \textstyle{\frac{(1-2k/n)}{n}
\biggl(  
3\binom{k+1}{2}+3\binom{k-n/3+1}{2}+3\sum_{j=2}^{\theb}
j(j+1)\binom{k+1-(\frac{1}{2}-\frac{1}{3j(j+1)})3n}{2} 
\biggr)
/n^2
+\Theta(n^3)} \\
&=24\tbinom{n}{4}\biggl(
\sum_{k=1}^{(n-2)/2}
\left(\frac{3}{2}\right)  \frac{\left(1-2k/n\right)}{n}  \left(  \frac{k}
{n}\right)  ^{2}
+
\sum_{k=n/3}^{(n-2)/2}
\left(  \frac{3}{2}\right)  \frac{\left(  1-2k/n\right)}{n} \left(  \frac{k}
{n}-\frac{1}{3}\right)^{2}\\
&{ } \ +
\sum_{j=2}^{\infty}
j(j+1)
\sum_{k=\left(
(1/2)-1/(3j(j+1))\right)n}^{(n-2)/2}
\left(  \frac{3}{2}\right)  
\frac{\left(  1-2k/n\right)}{n}
\left(
\frac{k}{n}
-
\left(
\frac{1}{2}-\frac{1}{3j(j+1)}
\right)
\right)^{2}
 \biggr)
+\Theta(n^{3}) \\
&
=24\tbinom{n}{4} \biggl(  \int_{0}^{1/2}(3/2)(1-2x)x^{2}dx+\int_{1/3}
^{1/2}(3/2)(1-2x)(x-1/3)^{2}dx \\
& +
{\textstyle\sum_{j=2}^{\infty}}
(3/2)j(j+1)\int_{\frac{1}{2}-\frac{1}{3j(j+1)}}^{1/2}
(1-2x)(x-(1/2-1/(3j(j+1)))^{2}dx\biggr)  +\Theta(n^{3})
\end{align*}

Elementary calculations show that $\frac{3}{2}\int_{0}^{1/2}(1-2x)x^{2}%
dx=3/8,\frac{3}{2}\int_{1/3}^{1/2}(1-2x)(x-1/3)^{2}dx=1/216,$ and $\int
_{\frac{1}{2}-\frac{1}{3j(j+1)}}^{1/2}(1-2x)(x-(1/2-1/(3j(j+1)))^{2}%
dx=(1/486)j^{4}(1+j)^{4}$.

Thus,
$$
\rcr(\dd_n) \ge 
24\tbinom{n}{4}\left(  3/8+1/216+(3/2)%
{\sum_{j=2}^{\infty}}
\frac{1}{486j^{3}(j+1)^{3}}\right)  +\Theta(n^{3}).
$$

Since 
${\sum_{j=2}^{\infty}}
\frac{1}{j^{3}(j+1)^{3}}=%
{\sum_{j=2}^{\infty}}
\left(  \frac{1}{j^{3}}-\frac{3}{j^{2}}+\frac{6}{j}-\frac{1}{(j+1)^{3}}%
-\frac{3}{(j+1)^{2}}-\frac{6}{j+1}\right)  =\frac{79}{8}-\pi^{2}$, then
$$
\rcr(\dd_n) \ge
\frac{2}{27}\biggl(15-\pi^{2}\biggr)\binom{n}{4}+\Theta(n^{3}). 
$$
\vglue -0.6 cm \rightline{\hal}

\section{Proof of 
Theorem~\ref{thm:mainksets}}\label{sec:proofmainksets}

The first step to prove Theorem~\ref{thm:mainksets} is to obtain an equivalent
(actually, more general) formulation in terms of circular
sequences (namely Proposition~\ref{prop:main} below). 

\subsection{Circular sequences: reducing Theorem~\ref{thm:mainksets}
  to Proposition~\ref{prop:main}}

All the geometrical information of a point set $S$ gets encoded in
(any halfperiod of) the {\em circular sequence} associated to $S$.  We
recall that a circular sequence on $n$ elements is a doubly infinite
sequence $\ldots \pi_{-1}, \pi_0, \pi_1, \ldots$ of permutations of
the points in $S$, where consecutive permutations differ in a
transposition of neighboring elements, and, for every $i$, $\pi_i$ is
the reverse permutation of $\pi_{i+{n\choose 2}}$. Thus a circular
sequence on $n$ elements has period $2{n\choose 2}$, and all the
information is encoded in an $n$--{\em halfperiod}, that is, a
sequence of  ${n\choose 2}+1$ consecutive permutations.

Each $n$--point set $S$ has an associated circular sequence ${\mathbf
  \Pi}_S$, which contains all the geometrical information of
  $S$~\cite{goodmanpollack}. As we observed above, any $n$--halfperiod
  $\Pi$ of ${\mathbf \Pi}_S$ contains all the information of ${\mathbf
  \Pi}_S$, and so $n$--halfperiods are usually the object of choice to
  work with.  In an $n$--halfperiod $\Pi = \pi_0, \pi_1, \ldots,
  \pi_{n\choose 2}$, the {\em initial} permutation is $\pi_0$   and  the {\em final permutation}
  is $\pi_{n\choose 2}$.

 Not every $n$--halfperiod $\Pi$ arises from a point set $S$.  We refer
 the reader to the seminal work by Goodman and
 Pollack~\cite{goodmanpollack} for further details.

Observe that if $S$ is $3$--decomposable, 
then there is 
an $n$--halfperiod $\Pi$ of the circular sequence associated to $S$,
whose points can be labeled $a_1, \ldots, a_{n/3}, b_1, \ldots,
b_{n/3}, c_1, \ldots, c_{n/3}$, so
that: 
\begin{description}
\item{(i)} The initial permutation $\pi_0$ reads 
$a_{n/3}, a_{n/3-1}, \ldots, a_1,
b_1, b_2, \ldots, b_{n/3}, c_1, c_2, \ldots, c_{n/3}$; 
\item{(ii)} there is an $s$
such that in the $(s+1)$--st permutation first the $b$'s appear
consecutively, then the $a$'s appear consecutively, and then the $c$'s
appear consecutively; and
\item{(iii)} there is a $t$, with $t > s$, 
such that in the $(t+1)$--st permutation first the $b$'s appear
consecutively, then the $c$'s appear consecutively, and then the $a$'s
appear consecutively.
\end{description}

\noindent{\bf Definition }
An $n$--halfperiod $\Pi$ that satisfies properties (i)--(iii) above is
$3$--{\em decomposable}.

A transposition that occurs between elements in sites $i$ and
$i+1$ is an $(i,i+1)$--{\em transposition}. An $i$--{\em critical}
tranposition is either an $(i,i+1)$--transposition or an $(n-i,
n-i+1)$--transposition, and a $(\le k)$--{\em critical} transposition
is a transposition that is $i$--critical for some $i \le k$.  If $\Pi$
is an $n$--halfperiod,  then $\Nomberle{k}{\Pi}$ denotes
the number of $(\le k)$--critical transpositions in $\Pi$. 

The key result is the following.

\begin{proposition}\label{prop:main}
Let $\Pi$ be a $3$--decomposable $n$--halfperiod, 
and let $k < n/2$. 
Then 
$$
\Nomberle{k}{\Pi} \ge 
\ygriega{k}{n}.
$$

\end{proposition}

\noindent{\bf Proof of Theorem~\ref{thm:mainksets}.}
Let $S$ be $3$--decomposable, and let $\Pi$ be an $n$--halfperiod of
the circular sequence associated to $S$,
that satisfies properties (i)--(iii) above. Then $\Pi$ is
$3$--decomposable. Now, for any point set $T$ and any halfperiod $\Pi_T$
associated to $T$, the $(\le k)$--critical transpositions of $\Pi_T$
are in one--to--one correspondence with $(\le k)$--sets of $T$.
Applying this to $\Pi$ and $S$, it follows that $\nomberle{k}{S} =
\Nomberle{k}{\Pi}$.  Applying Proposition~\ref{prop:main},
Theorem~\ref{thm:mainksets} follows. \hal

\bigskip

We devote the rest of this section to the proof of
Proposition~\ref{prop:main}.

\subsection{Proof of Proposition~\ref{prop:main}}

Throughout this section, $\Pi=(\pi_0, \pi_1, \ldots, \pi_{n\choose
  2})$ is a $3$--decomposable $n$--halfperiod, with initial
  permutation
$\pi_0 = (a_{n/3}, a_{{n/3}-1}, \ldots, a_1, b_1, \ldots, b_{n/3}, c_1, \ldots,
  c_{n/3})$.   

In order to (lower) bound the number of $(\le k)$--critical
transpositions in $3$--decomposable circular sequences, we distinguish
between two types of transpositions. A transposition is {\em
  homogeneous} if it occurs between two $a$'s, between two $b$'s, or
between two $c$'s; otherwise it is {\em heterogeneous}. 
We let $\Nomberlehom{k}{\Pi}$ (respectively $\Nomberlehet{k}{\Pi}$)
denote the number of homogeneous (respectively heterogeneous) $(\le
k)$--critical transpositions in $\Pi$, so that
\begin{equation}\label{eq:suma}
\Nomberle{k}{\Pi} = \Nomberlehom{k}{\Pi} + \Nomberlehet{k}{\Pi}.
\end{equation}

\subsubsection{Bounding (actually, calculating) $\Nomberlehet{k}{\Pi}$}

Let us call a transposition an $ab$--{\em transposition} if it
involves one $a$ and one $b$. We similarly define $ac$-- and
$bc$--transpositions.  Thus, each heterogeneous transposition is
either an $ab$-- or an $ac$-- or a $bc$--transposition.

Since in $\Pi$ each $ab$--transposition moves the involved $a$ to the
right and the involved $b$ to the left, 
then (a) for each $i \le n/3$, there are {\em exactly} $i$
$i$--critical $ab$ transpositions; and (b) for each $i$, $n/3 < i <
2n/3$, there are {\em exactly} $2n/3-i$ $i$--critical
$ab$--transpositions.  Since the same holds for $ac$-- and
$bc$--transpositions, it follows that for each $i \le n/3$, there are
{\em exactly} $3i$ $i$--critical heterogeneous transpositions, and for
each $i$, $n/3<i<2n/3$, exactly $2n-3i$ $i$--critical heterogeneous
transpositions.

\def\tfl#1{{\lfloor{#1}\rfloor}}
\def\tcl#1{{\lceil{#1}\rceil}}

Therefore, for each $k \le n/3$, there are exactly 
$\sum_{i\le k} 3i 
= 3\binom{k+1}{2}$ $(\le k)$--critical transpositions, and  
if $n/3 < k < n/2$, then there are exactly 
$\sum_{i\le n/3} 3i 
+
\sum_{n/3 < i \le k} 2n-3i
+
\sum_{n-k-1 < i \le 2n/3 -1} 2n-3i
=
3\binom{n/3+1}{2}
+
(k-n/3)n$
$(\le k)$--critical transpositions.

We now summarize these results.

\begin{proposition}\label{pro:heterogeneous}
Let  $\Pi$ be a $3$--decomposable $n$--halfperiod, and let $k < n/2$. 
Then 
\begin{equation*}
\Nomberlehet{k}{\Pi} = 
\begin{cases} 3\binom{k+1}{2} & \hbox{\hglue 0.3 cm} \text{if $k \le n/3$,}
\\[0.4 cm]
3\binom{n/3+1}{2} + (k-n/3)n 
& \hbox{\hglue 0.3 cm} \text{if $n/3 < k < n/2$,} 
\end{cases}
\end{equation*}
\end{proposition}

\subsubsection{Bounding $\Nomberlehom{k}{\Pi}$}

Our goal here is  to give a lower bound 
(see Proposition~\ref{pro:homogeneous}) for the number 
$\Nomberlehom{k}{\Pi}$ of
homogeneous $(\le k)$--critical transpositions in a $3$--decomposable $n$--halfperiod
$\Pi$. 

Our approach is to find an {\em upper} bound for 
$\Momberaa{k}{\Pi}$, which will denote the number of
$aa$--transpositions that are {\em not} $(\le k)$--critical
($\Momberbb{k}{\Pi}$ and $\Mombercc{k}{\Pi}$ are defined analogously).  Since
the total number of $aa$--transpositions is $\binom{n/3}{2}$, this
will yield a lower bound for the contribution of $aa$--transpositions
(and, by symmetry, for the contribution of $bb$--transpositions and of
$cc$--transpositions) 
to $\Nomberlehom{k}{\Pi}$.

\begin{remark}\label{rem:couldbe0}
For every $k \le n/3$, it is a trivial task to construct $n$--halfperiods
$\Pi$ for which $\Nomberlehom{k}{\Pi} = 0$.  
In view of this, we concentrate our efforts
on the case $k > n/3$. 
\end{remark}

A transposition between elements in positions $i$ and $i+1$, with $k+1
\le i \le n-k-1$, is {\em valid}.  Thus our goal is to (upper) bound the
number of valid $aa$--transpositions.

Let $\digaawithpi$ be the digraph with vertex set $a_1, \ldots, a_{n/3}$, and
such that there is a directed edge from $a_\ell$ to $a_j$ if and only
if $\ell < j$ and the transposition that swaps $a_\ell$ and $a_j$ is
valid.  For $j=1,\ldots, n/3$, we let $\outdpi{a_j}$ (respectively
$\indpi{a_j}$) denote the outdegree (respectively indegree) of $a_j$ in
$\digaawithpi$. 
We define $\digabwithpi, \digacwithpi, \indpi{b_j}, \outdpi{b_j}, \indpi{c_j}$ 
and $\outdpi{c_j}$ analogously. 

The inclusion of the symbol $\Pi$ in $\digaawithpi, \indpi{a_i}$,
etc., is meant to emphasize the dependence on the specific
$n$--halfperiod $\Pi$. For brevity we will omit the reference to $\Pi$
and simply write $\digraa, \digab, \digac, \ind{a_i}, \outd{a_i}$,
and so on. No confusion will arise from this practice.

The importance of $\digraa, \digab$, and $\digac$ 
is clear from the following observation.

\begin{remark}\label{rem:digraph}
For each $n$--halfperiod $\Pi$, the number of edges of $\digraa$ {\em equals}
$\Momberaa{k}{\Pi}$.  Indeed, to each valid $aa$--transposition, that is,
each transposition that contributes to $\Momberaa{k}{\Pi}$, there
corresponds a unique edge in $\digraa$.  Analogous observations hold
for $\digab$ and $\digac$.
\end{remark}

In view of Remark~\ref{rem:digraph}, we direct our efforts to bounding
the number of edges in $\digraa$.  The essential observation to get
this bound is the following:

\begin{equation}\label{eq:degrees}
\ind{a_j} \le \min\{ 
n-2k-1 + \outd{a_j}
,
(n/3)-j
\}. 
\end{equation} 

To see this, simply note that , $\ind{a_j} \le n-2k-1 + \outd{a_j}$,
since $n-2k-1 + \outd{a_j}$ is clearly the maximum possible number of
valid moves in which $a_j$ moves right, and trivially $\ind{a_j} \le
(n/3)-j$, since there are only $(n/3)-j$ $a_\ell$'s with $\ell < j$.

\begin{proposition}\label{pro:boundingnoedges}
If $\Pi$ is a $3$--decomposable $n$--halfperiod, and $n/3 < k < n/2$,
then $\digraa$ has at most $\binom{n/3}{2} - (1/3)\left(\ygriega{k}{n}
- 3\binom{n/3+1}{2} - (k-n/3)n\right)$edges.
\end{proposition}

\noindent{\em Proof.}  
Let $\dd_{k,n}$ denote the class of all
digraphs
with vertex set
$a_1, \ldots, a_{n/3}$, with every directed edge $\ora{a_\ell a_j}$ 
satisfying $\ell < j$ and 
$\ind{a_j} \le \min\{ 
n-2k-1 + \outd{a_j}
,
{n/3}-j
\}$.

We argue that any graph in $\dd_{k,n}$ has at most $\binom{n/3}{2} -
(1/3)\bigl(\ygriega{k}{n}$ $ - 3\binom{n/3+1}{2} - (k-n/3)n\bigr)$
edges.  This clearly finishes the proof, since $\digraa \in \dd_{k,n}$.

To achieve this, we note that it
follows from the work in Section 2
in~\cite{baloghsalazar} that the maximum number of edges of such a
digraph is attained in the digraph $D_{k,n}$ recursively constructed
as follows.
First define that all the directed
edges arriving at $a_{n/3}$ are the edges $\ora{a_j a_{n/3}}$ for
$j=(n/3)-1,\ldots,(n/3)-n-2k-2$.  Now, for $j +1 \le n/3$, once all the directed
edges arriving at $a_{j+1}$ have been determined, fix that (all) the
directed edges arriving at $a_j$ are $\ora{a_\ell a_j}$, for all those
$\ell$ that satisfy $j - \ell \le  n-2k-1 + \outd{a_j}$.

Since no digraph in $\dd_{k,n}$ has more edges than $D_{k,n}$, to
finish the proof it suffices to bound the number of edges of
$D_{k,n}$. This is the content of Claim~\ref{cla:theclaim} below.

\medskip

\begin{claim}\label{cla:theclaim}
$D_{k,n}$ has at most 
$
\binom{n/3}{2} - 
\bigl(
\ygriega{k}{n}$ $ - 3\binom{n/3+1}{2} - (k-n/3)n
\bigr)/3$
edges.
\hfill 
\end{claim}

\medskip

\noindent{\em Sketch of proof.}  Since we know the exact indegree of
each vertex in $D_{k,n}$, we know the exact number of edges of
$D_{k,n}$, and so the proof of Claim~\ref{cla:theclaim} is no more
than a straightforward, but quite long and tedious, calculation. \hal

\bigskip

\begin{corollary}\label{cor:alsobandc}
If $\Pi$ is a $3$--decomposable $n$--halfperiod, and $n/3 < k < n/2$,
then each of $\digab$ and $\digac$ has at most $\binom{n/3}{2} -
(1/3)\left(\ygriega{k}{n} - 3\binom{n/3+1}{2} - (k-n/3)n\right)$edges.
\end{corollary}

\noindent{\em Proof.} In the proof of
Proposition~\ref{pro:boundingnoedges}, the only
 relevant property about $\digraa$ is that the a's
form a set of $n/3$ points that in some permutation of
$\Pi$ (namely $\pi_0$)
 appear all consecutively and at
the beginning of the permutation. Since $\Pi$ is $3$--decomposable,
this condition is also satisfied by the set of b's and by the set of  c's.  \hal

\bigskip

We now summarize the results in the current subsection.

\begin{proposition}\label{pro:homogeneous}
If $\Pi$ is a $3$--decomposable $n$--halfperiod, and 
$n/3 < k < n/2$, then
$$
\Nomberlehom{k}{\Pi}  \ \ge \ 
\ygriega{k}{n} - 3\binom{n/3+1}{2} - (k-n/3)n.$$
\end{proposition}

\noindent{\em Proof.} By Remark~\ref{rem:digraph}, the number
$\Momberaa{k}{\Pi}$ of
$aa$--transpositions that are {\em not} 
$(\le k)$--critical equals the number of
edges in $\digraa$, which by Proposition~\ref{pro:boundingnoedges} is at most
$(1/3)(\ygriega{k}{n} - 3\binom{n/3+1}{2} - (k-n/3)n)$.
Since the total number of $aa$--transpo\-sitions
is $\binom{n/3}{2}$, then  the number of
$aa$--transpositions that contribute to $\Nomberle{k}{\Pi}$ is at
least $\binom{n/3}{2} - 
\bigl(\binom{n/3}{2} - 
(1/3)\bigl(\ygriega{k}{n} - 
3\binom{n/3+1}{2} - (k-n/3)n\bigr)\bigr)
$
$ = (1/3)\left(\ygriega{k}{n} - 
3\binom{n/3+1}{2} - (k-n/3)n\right)$.  
By Corollary~\ref{cor:alsobandc}, $bb$-- and 
$cc$--transpositions contribute in at least the same amount
to
$\Nomberlehom{k}{\Pi}$, and so the claimed inequality follows. \hal

\subsubsection{Proof of Proposition~\ref{prop:main}}

Proposition~\ref{prop:main} follows immediately from
Eq.~(\ref{eq:suma}) and Propositions~\ref{pro:heterogeneous}
and~\ref{pro:homogeneous}.  \hal

\section{Concluding remarks}\label{sec:concludingremarks}

All the lower bounds proved above remain true for point sets that
satisfy conditions (i) and (ii) (and not necessarily condition (iii))
for $3$--decomposability.

  \renewcommand{\theequation}{A-\arabic{equation}}
  \setcounter{equation}{0}  

\bigskip

\newpage

\noindent{\bf \Large Appendix: Proof of Claim~\ref{cla:theclaim}}

\bigskip
\medskip

\setcounter{section}{0}

\setcounter{theorem}{0}


Since $D_{k,n}$ is a well--defined digraph, and we know the exact
indegree of each of its vertices, Claim~\ref{cla:theclaim} is no more
than long and tedious, yet  elementary, calculation.

The purpose of this Appendix is to give the full details of this
calculation.

We prove Claim~\ref{cla:theclaim} in two steps. First we obtain an
expression for the {\em exact} value of the number of edges of
$D_{k,n}$, and then we show that this exact value is upper bounded by
the expression in Claim~\ref{cla:theclaim}.

\section{The exact number of edges in $D_{k,n}$}

The exact number of edges in $D_{k,n}$ is a function of the following
parameters. Let $i,j$ be positive integers with $i \le j$. Then:
\begin{itemize}
\item $\theb(i,j)$ is the (unique) nonnegative integer such that
$\tbinom{\theb(i,j)+1}{2} <j/i \leq \tbinom{\theb(i,j)+2}{2}$; and
\item $\theq(i,j)$ and $\ther(i,j)$ are the (unique) integers satisfying
$0\leq \theq(i,j) < i,1\leq \ther(i,j)\leq \theb(i,j)+1$ and such that
\begin{equation}\label{eq:a1}
j=i \tbinom{\theb(i,j)+1}{2}+\theq(i,j)(\theb(i,j)+1)+\ther(i,j)  
\end{equation}
\end{itemize}

For brevity, in the rest of the section we let $\enetercios:=n/3$ and
$m:=n-2k-1$.

The key observation is that we know the indegree of each vertex in
$D_{k,n}$:

\begin{proposition}[Proposition 19 in~\cite{baloghsalazar}] \label{ingrade}
For each integer $1\le i \le \enetercios$, and each vertex $a_i$ of $D_{k,n}$,
$\ind{a_i}=\theb(i,\enetercios) m +\theq(i,\enetercios)$.
\end{proposition}

The number of edges of $D_{k,n}$ equals the sum of the
indegrees over all vertices in $D_{k,n}$.   Thus our main task is to
find a closed expression for the sum 
$\sum_{1\leq i \leq \enetercios} \ind{a_i}$.  
This is the content of
our next statement.

\begin{proposition}[Exact number of edges of $D_{k,n}$]\label{numberofedges}
The number $\sum_{1 \le i \le s} \ind{a_i}$ of edges of $D_{k,n}$ is
\begin{align}
E(k,n) \hbox{\hglue 0.1 cm}  :=  \  & \hbox{\hglue 0.2 cm} 
2m^2
\binom{\theb(m,\enetercios)+1}{3}+ \binom{\theb(m,\enetercios)+1}{2}
\binom{m}{2}+ 2m\cdot\theq(m,\enetercios)
\binom{\theb(m,\enetercios)+1}{2}+ \nonumber \\
&
\binom{\theq(m,\enetercios)}{2}\binom{\theb(m,\enetercios)+1}{1}+\ther(m,\enetercios)
\biggl(m
\cdot\theb(m,\enetercios)+\theq(m,\enetercios)\biggr). \nonumber
\end{align}
\end{proposition}

\noindent{\em Proof.}  
We break the index set of the summation 
$\sum_{1\leq i \leq \enetercios} \ind{a_i}$ into three parts, in terms
of $\alpha:=m\tbinom{\theb(m,\enetercios)+1}{2}$ and 
$\beta:=\theq(m,\enetercios) (\theb(m,\enetercios)+1)$.  
We let  $A= \sum_{1\leq
i\leq\alpha}\ind{a_i}$, $B= \sum_{\alpha+1\leq
i\leq\alpha+\beta}\ind{a_i}$, and $C=\sum_{\alpha+\beta+1\leq i\leq
\enetercios}\ind{a_i}$ so that 
\begin{equation}\label{eq:decomp}
\sum_{1\leq i\leq\enetercios}\ind{a_i}=A+B+C.
\end{equation}

We calculate each of $A$, $B$, and $C$
separately.

\bigskip
\noindent{\em Calculating $A$}
\bigskip

If $\ell,j$ are integers such that $0\leq \ell \leq m-1$
and $0\leq j\leq \theb(m,\enetercios)$, we define $S_{j}:=\{i:\theb(i,m)=j\}$
and $T_{j,\ell} :=\{i:\theb(i,m)=j,\theq(i,m)=\ell\}.$ Note that
$S_{1},S_{2},...,S_{\theb(m,\enetercios)}\ $is a partition of
$\{1,2,...,n\}$ and that for each $j \leq
\theb(m,\enetercios)-1,T_{j,0},T_{j,1},...,T_{j,m-1}$ is a partition of
$S_{j}.$

Note that $A$ can be rewritten as $
\sum_{0\leq j\leq \theb(m,\enetercios)-1}\sum_{i\in S_{j}}\ind{a_i}.$ By
Proposition \ref{ingrade} this equals $\sum_{0\leq j\leq
\theb(m,\enetercios)-1} \sum_{i\in S_{j}}(m\cdot\theb(i,m)+\theq(i,m))$. That is,

$$ A= \sum_{0\leq j\leq \theb(m,\enetercios)-1} \left( m \sum_{i\in S_{j}}
\theb(i,m)+\sum_{i\in S_{j}} \theq(i,m)\right). $$

Since $0\leq \theq(i,m) \leq m-1$ for all $i$, and $T_{j,0},$
$T_{j,1},...,T_{j,m-1}$ is a partition of $S_{j}$, then $\sum _{i\in
S_{j}} \theq(i,m)= \sum _{0 \leq l \leq m-1} \sum _{i\in T_{j,l}} \theq(i,m).$
Thus,

\begin{equation} \label{A}
A= \sum _{0\leq j\leq \theb(1,\enetercios)-1}
\left(  m \sum _{i \in S_{j}}  \theb(i,m)+ \sum _{0\leq l\leq m-1}
\sum _{i\in T_{j,l}} \theq(i,m) \right). 
\end{equation}

On other hand, for $0\leq j\leq \theb(m,\enetercios)-1$ and $0\leq \ell\leq
m-1$, it is not difficult to verify that $\left\vert T_{j,\ell}\right\vert
$ $=j+1$. This implies that $\left\vert S_{j}\right\vert =m(j+1).$

By definition of $S_{j}$ we have 

\begin{equation} \label{S_j}
\sum _{i \in S_{j}} \theb(i,m)= \sum _{i\in S_{j}} j = 
 j \left\vert S_{j} \right\vert =jm(j+1).
\end{equation}

By definition of $T_{j,\ell}$ we have

\begin{equation}\label{T_jl}
\sum _{i\in T_{j,\ell}} \theq(i,m)= \sum _{i\in T_{j,\ell}}\ell= 
\ell\left\vert T_{j,\ell}\right\vert =\ell(j+1).
\end{equation}

Substituting (\ref{S_j}) and (\ref{T_jl}) into (\ref{A}) we obtain
\begin{align}
A \  & = {\sum _{0\leq j\leq \theb(m,\enetercios)-1}} \left( m(jm(j+1))+
{\sum _{0\leq \ell\leq m-1}} \ell (j+1) \right) \nonumber \\
& = {\sum _{0\leq j\leq \theb(m,\enetercios)-1}} \left(
2m^{2}\tbinom{j+1}{2}+(j+1) {\sum _{0\leq \ell\leq m-1}} \ell\right) \nonumber \\
& ={\sum _{0\leq j\leq \theb(m,\enetercios)-1}} \left(
2m^{2}\tbinom{j+1}{2}+(j+1)\tbinom{m}{2}\right) \nonumber \\[0.2 cm]
& =2m^{2}\tbinom{\theb(m,\enetercios)+1}{3}+
\tbinom{\theb(m,\enetercios)+1}{2}\tbinom{m}{2}. \label{eq:fora}
\end{align}

\bigskip
\noindent{\em Calculating $B$}
\bigskip

Since $\theb(i,m)=\theb(m,\enetercios)$ for each $i \geq \alpha+1$,
and $\ind{a_i}=$ $m\cdot\theb(i,m)+\theq(i,m)$, then
$B= {\sum _{\alpha+1\leq i\leq\alpha+\beta}} \ind{a_i}=
{\sum _{\alpha+1\leq i\leq\alpha+\beta}}( m\cdot\theb(i,m)+\theq(i,m)).$ 

\bigskip
Therefore
\begin{align}
B & = {\sum _{\alpha+1\leq i\leq\alpha+\beta}} m\cdot\theb(m,\enetercios)+
{\sum_{\alpha+1\leq i\leq\alpha+\beta}} \theq(i,m) \nonumber \\
& =m\cdot\theb(m,\enetercios) {\sum _{\alpha+1\leq i\leq\alpha+\beta}}1+ {\sum
 _{\alpha+1\leq i\leq\alpha+\beta}} \theq(i,m) \nonumber   \\
& =m\cdot\theb(m,\enetercios)\beta+{\sum _{\alpha+1\leq i\leq\alpha+\beta}}
\theq(i,m) \nonumber  \\
&=m\cdot\theb(m,\enetercios)\theq(m,\enetercios)(\theb(m,\enetercios)+1)+ {\sum
_{\alpha+1\leq i\leq\alpha+\beta}} \theq(i,m). \nonumber 
\end{align}

On other hand it is easy to check that $\left\vert
T_{\theb(m,\enetercios) ,k}\right\vert =\theb(m,\enetercios)+1$ for every
 $k$ such that $0\leq k\leq \theq(m,\enetercios)-1.$ Since $0\leq
\theq(i,m)\leq \theq(m,\enetercios)-1$ for every $i$ such that $\alpha+1\leq
i\leq\alpha+\beta,$ then $T_{\theb(m,\enetercios),0},$
$T_{\theb(m,\enetercios),1} ,...,T_{\theb(m,\enetercios),\theq(m,\enetercios)-1}$
is a partition of $\{\alpha+1,\alpha +2,...,\alpha+\beta\}.$ Thus,

\begin{equation}\label{B}
{\sum _{\alpha+1\leq i\leq\alpha+\beta}} \theq(m,\enetercios)= {\sum
_{0\leq \ell\leq \theq(m,\enetercios)-1}} \hbox{\hglue 0.2 cm} {\sum _{i\in
T_{\theb(m,\enetercios),\ell}}} \theq(i,m).
\end{equation}

We note that $ {\sum _{i\in T_{\theb(m,\enetercios),\ell}}} \theq(i,m)=
\ell\left\vert T_{\theb(m,\enetercios),k}\right\vert =\ell(\theb{(m,\enetercios)}+1).$
Using this fact  in (\ref{B}) we obtain

$${\sum _{\alpha+1\leq i\leq\alpha+\beta}}\theq(i,m)={\sum _{0\leq \ell\leq
\theq{(m,\enetercios)}-1}}
\ell(\theb{(m,\enetercios)}+1)=\tbinom{\theq{(m,\enetercios)}}{2}(\theb{(m,\enetercios)}+1).$$

Thus,
\begin{align}
B & =m\cdot\theb{(m,\enetercios)}\theq{(m,\enetercios)}(\theb(m,\enetercios)+1)+
\tbinom{\theq{(m,\enetercios)
}}{2}(\theb{(m,\enetercios)}+1) \nonumber \\ 
&=2m\cdot\theq{(m,\enetercios)}\tbinom{\theb{(m,\enetercios)}+1}{2}
+\tbinom{\theq{(m,\enetercios)}}{2}(\theb{(m,\enetercios)}+1). \label{eq:forb}
\end{align}

\bigskip
\noindent{\em Calculating $C$}
\bigskip

Since $\theb(i,m)=\theb(m,\enetercios)$;
$\theq(i,m)=\theq(m,\enetercios)$ for each  $i$ such that
$i\geq\alpha+\beta+1$; and $\ind{a_i}=$ $m\cdot\theb(i,m)+\theq(i,m)$, it follows that
\begin{align}
 C=\sum_{\alpha+\beta+1\leq i\leq \enetercios} \ind{a_i}  
& = \hbox{\hglue 0.25 cm}  \sum_{\alpha+\beta+1\leq i\leq \enetercios} m\cdot\theb(i,m)+\theq(i,m)
 \nonumber \\
&=  \hbox{\hglue 0.25 cm}   {\sum_{\alpha+\beta+1\leq i\leq \enetercios}}
m\cdot\theb(m,\enetercios)+\theq(m,\enetercios) \nonumber \\
&= \hbox{\hglue 0.25 cm}  (\enetercios-\alpha-\beta)
\biggl(m\cdot\theb(m,\enetercios)+\theq(m,\enetercios)\biggr)
\nonumber
\end{align}

From (\ref{eq:a1}) it follows that
$\ther(m,\enetercios)=(\enetercios-\alpha-\beta)$, and so
\begin{equation}\label{eq:forc}
C=\ther(m,\enetercios)(m\cdot\theb(m,\enetercios)+\theq(m,\enetercios)).
\end{equation}

Now from (\ref{eq:fora}), (\ref{eq:forb}), and (\ref{eq:forc}), it
 follows that $E(k,n) = A + B + C$, and 
so Proposition~\ref{numberofedges} follows from
 (\ref{eq:decomp}). \hal

\section{Upper bound for number of edges in $D_{k,n}$: \\ Proof 
of Claim~\ref{cla:theclaim}}

First we bound the number of $(\le k)$--edges in $3$--decomposable
$n$--halfperiods in terms of the expression $E(k,n)$ in
Proposition~\ref{numberofedges}.

\begin{proposition}\label{exact}
Let $\Pi$ be a $3$--decomposable $n$--halfperiod, and let $k <
n/2$. Then
\begin{equation*}
\Nomberle{k}{\Pi} \geq  
\begin{cases} 3\binom{k+1}{2} & \hbox{\hglue 0.3 cm} \text{if $k \le n/3$,}
\\[0.4 cm]
3\binom{n/3+1}{2} + (k-n/3)n+3 \biggl( \binom{n/3}{2} - E(k,n)\biggr)
& \hbox{\hglue 0.3 cm} \text{if $n/3 < k < n/2$.} \\[0.4 cm]
\end{cases}
\end{equation*}

\end{proposition}

\noindent{\em Proof.} Obviously, $\Nomberle{k}{\Pi} \geq
\Nomberlehet{k}{\Pi}$ and so the case $k \le n/3$ follows from
Proposition~\ref{pro:heterogeneous}.  Now suppose that $n/3 < k <
n/2$.  Recall that $\Nomberlehom{k}{\Pi} = \Momberaa{k}{\Pi} +
\Momberbb{k}{\Pi} + \Mombercc{k}{\Pi}$. Now the total number of
$aa$-- (and $bb$--, and $cc$--) transpositions is exactly
$\binom{n/3}{2}$, and so $\Nomberlehom{k}{\Pi} =3 \binom{n/3}{2} -
\Momberaa{k}{\Pi} - \Momberbb{k}{\Pi} - \Mombercc{k}{\Pi}$. 
Thus it follows from Remark~\ref{rem:digraph} and
Proposition \ref{numberofedges} that $\Nomberlehom{k}{\Pi} \geq 3(
\binom{n/3}{2} - E(k,n)$. This fact, together\marginpar{\tiny Must
  explain that for $bb$ and $cc$ same digraph works} with Proposition
\ref{pro:heterogeneous}, implies that $\Nomberle{k}{\Pi}=
\Nomberlehet{k}{\Pi}+ \Nomberlehom{k}{\Pi}\geq 3\binom{n/3+1}{2} +
(k-n/3)n+3( \tbinom{n/3}{2} - E(k,n))$, as claimed.  \hal

\bigskip
\bigskip

\noindent{\bf Proof of Claim~\ref{cla:theclaim}.}  Recall 
that $\enetercios:=n/3$ and
$m:=n-2k-1$. By  Remark~\ref{rem:couldbe0}  we know that
 $k>n/3$,  and so $\enetercios \geq m$.  
From (\ref{eq:a1}), it follows that
\begin{equation}\label{q}
\theq(m,\enetercios)=\frac{\enetercios-m \tbinom{\theb(m,\enetercios)+1}{2}-
\ther(m,\enetercios)}{\theb(m,\enetercios)+1}.
\end{equation}

Now by Proposition \ref{exact}, $\Nomberle{k}{\Pi} \geq L(k,n), $ where
\begin{equation}\label{L}
L(k,n):=3\binom{\enetercios+1}{2} + (k-\enetercios)n+3(
\tbinom{\enetercios}{2} - E(k,n)).
\end{equation}

Substituting in  $E(k,n)$ the value of $\theq(m,\enetercios)$ 
given in (\ref{q}), a (long and tedious yet) totally elementary simplification
yields 
\begin{align}
\nonumber
L(n,k)- \ygriega{k}{n}-1/3 
 = & \frac{1}{8( \theb(m,\enetercios)+1)} \biggl(
5\theb(m,\enetercios)^2 + 4\theb(m,\enetercios)^3 + \theb(m,\enetercios)^4 \\
&  
+ \theb(m,\enetercios) ( - 12\ther(m,\enetercios)+2) + 
\nonumber 
12(\ther(m,\enetercios)-1)\ther(m,\enetercios)\biggr). 
\end{align}

Define $f(\theb(m,\enetercios), \ther(m,\enetercios)):=
L(n,k)- \ygriega{k}{n}-1/3$.
An elementary calculation shows that $\ther_0(m,\enetercios):=(\theb(m,\enetercios)+1)/2$
minimizes $f(\theb(m,\enetercios), \ther(m,\enetercios))$. Thus
$f(\theb(m,\enetercios), \ther(m,\enetercios)) \geq f(\theb(m,\enetercios),
\ther_0(m,\enetercios))=
(\theb(m,\enetercios)+3)(\theb(m,\enetercios)+1)(\theb(m,\enetercios)-1)/8$. Since
$\theb(m,\enetercios)$ is a nonnegative integer, it follows that 
$f(\theb(m,\enetercios), \ther(m,\enetercios))\geq -1/3$ and therefore
$L(n,k)- \ygriega{k}{n}\geq 0$.

By  (\ref{L}), $E(k,n)= (1/3)(3 \binom{\enetercios+1}{2} +
(k-\enetercios)n+ 3 \tbinom{\enetercios}{2} - L(k,n)).$ Since $-L(n,k)
\leq -\ygriega{k}{n}$, then $E(k,n) \leq (1/3)(3
\binom{\enetercios+1}{2} + (k-\enetercios)n+ 3
\tbinom{\enetercios}{2}- \ygriega{k}{n})=
\tbinom{\enetercios}{2}-(1/3)(\ygriega{k}{n}- 3
\binom{\enetercios+1}{2}- (k-\enetercios)n ).$  This proves
Claim~\ref{cla:theclaim}, since $E(k,n)$ is the total number of 
edges in $\dd_{k,n}$. \hal


\begin{thebibliography}{99}

\bibitem{af} B.M.~\'Abrego and S.~Fern\'andez--Merchant, A lower bound
  for the rectilinear crossing number, {\it Graphs and Comb.},
{\bf 21} (2005), 293--300.

\bibitem{seq} B.M.~\'Abrego, J.~Balogh, S.~Fern\'andez--Merchant,
  J.~Lea\~nos, and G.~Salazar, An extended lower bound on the number
  of $(\le k)$-edges to generalized configurations of points and the
  pseudolinear crossing number of $K_n$. Submitted (2007).

\bibitem{agor} O.~Aichholzer, J.~Garc\'\i a, D.~Orden, and P.~Ramos,
  New lower bounds for the number of $({\le}{k})$-edges and the
  rectilinear crossing number of $K_n$, {\em Discr.~Comput.~Geom.}, to appear.

\bibitem{aweb} O.~Aichholzer. On the rectilinear crossing
  number. Available online at {\tt
http://www.ist.tugraz.at/ staff/aichholzer/crossings.html.}


\bibitem{baloghsalazar} J.~Balogh and G.~Salazar, $k$--sets, convex
  quadrilaterals, and the rectilinear crossing number of $K_n$, {\it
  Discr.~Comput.~Geom.} {\bf 35} (2006), 671--690.

\bibitem{bmp} P.~Brass, W.O.J.~Moser, and J.~Pach, {\it Research Problems
  in Discrete Geometry.} Springer, New York (2005).

\bibitem{goodmanpollack} J.~E. Goodman and R.~Pollack,  On the 
combinatorial classification of nondegenerate configurations in  
the plane,  {\it J. Combin. Theory Ser. A}  {\bf 29}  (1980), 220--235.


\bibitem{lovasz} L.~Lov\'asz, K.~Vesztergombi, U.~Wagner, and E.~Welzl, 
Convex Quadrilaterals and $k$--Sets. 
{\it Towards a Theory of Geometric Graphs}, (J.~Pach, ed.),
Contemporary Math., AMS, 
139--148 (2004).


\end{thebibliography}
\end{document}